\documentclass[11pt,oneside]{amsproc}
\usepackage[T2A]{fontenc}
\usepackage[russian]{babel}
\usepackage{amssymb}
\newtheorem{tm}{Теорема}
\newtheorem{lem}{Лемма}
\newtheorem{dif}{Определение}
\newtheorem{prop}{Утверждение}
\title[Оценки числа собственных значений оператор-функций]{Оценки числа
собственных значений абстрактных и дифференциальных оператор-функций}
\author{А.~А.~Владимиров}
\thanks{Работа поддержана РФФИ, гранты \No\,01-01-00691
и \No\,02-01-06464.}
\begin{document}
\begin{abstract}
В статье рассматривается оператор-функция \(F(\lambda)\),
\(\lambda\in(\sigma,\tau)\subseteq\mathbb R\), значениями которой являются
полуограниченные самосопряженные операторы пространства \(\mathfrak H\).
Для числа \(\mathcal N_F(\alpha,\beta)\) собственных значений такой
оператор-функции на полуинтервале \([\alpha,\beta)\Subset(\sigma,\tau)\)
устанавливаются оценки вида \(\mathcal N_F(\alpha,\beta)\geqslant
\nu_F(\beta)-\nu_F(\alpha)\) и равенства вида \(\mathcal
N_F(\alpha,\beta)=\nu_F(\beta)-\nu_F(\alpha)\), где \(\nu(\xi)\) есть число
отрицательных собственных значений оператора \(F(\xi)\),
\(\xi\in(\sigma,\tau)\).

Полученные абстрактные результаты применяются к обыкновенным дифференциальным
оператор-функциям на конечном интервале.
\end{abstract}
\maketitle

\section{Абстрактные утверждения}
Пусть \(\mathfrak H\) есть некоторое абстрактное гильбертово
пространство, и пусть \(F(\lambda)\) есть абстрактная оператор-функция
интервала \((\sigma,\tau)\subseteq\mathbb R\) со следующими свойствами.
\begin{enumerate}
\item При любом фиксированном \(\lambda_0\in(\sigma,\tau)\) оператор
\(F(\lambda_0)\) представляет собой действующий в пространстве \(\mathfrak H\)
самосопряженный оператор с компактной резольвентой. 
\item Для любого фиксированного \(\lambda_0\in(\sigma,\tau)\) найдутся такое
вещественное число \(\mu\) и такая окрестность \(\mathcal U(\lambda_0)\Subset
(\sigma,\tau)\) числа \(\lambda_0\), что при любом фиксированном \(\lambda_1
\in\mathcal U(\lambda_0)\) оператор \(F(\lambda_1)-\mu\) будет неотрицателен.
\item Оператор-функция \((F(\lambda)-i)^{-1}\) непрерывна на интервале
\((\sigma,\tau)\) в смысле равномерной операторной топологии.
\end{enumerate}

Как обычно, \emph{собственным значением кратности \(m\)} оператор-функции
\(F(\lambda)\) будет называться такое число \(\lambda_0\in(\sigma,\tau)\), для
которого \(0\) является собственным значением кратности \(m\) оператора
\(F(\lambda_0)\). Кроме того, через \(\mathcal N_F(\xi_1,\xi_2)\) будет
обозначаться сосчитанное с учетом кратности число собственных значений
оператор-функции \(F(\lambda)\), лежащих на полуинтервале
\([\xi_1,\xi_2)\Subset(\sigma,\tau)\).

Обозначим также через \(\nu_F(\lambda)\) числовую функцию, для которой при
любом фиксированном \(\lambda_0\in(\sigma,\tau)\) число \(\nu_F(\lambda_0)\)
представляет собой сосчитанное с учетом кратности число отрицательных
собственных значений оператора \(F(\lambda_0)\). В настоящем пункте будет
исследован вопрос о связи между поведением функции \(\nu_F(\lambda)\) и
поведением величин \(\mathcal N_F(\xi_1,\xi_2)\).

\subsection{Простейшая оценка числа собственных значений}
Справедливо следующее утверждение.
\begin{tm}\label{tm:1.44}
Для любого полуинтервала \([\xi_1,\xi_2)\Subset(\sigma,\tau)\) имеет место
оценка
\[
	\mathcal N_F(\xi_1,\xi_2)\geqslant\nu_F(\xi_2)-\nu_F(\xi_1).
\]
\end{tm}

Чтобы доказать теорему~\ref{tm:1.44}, введем в рассмотрение следующую
вспомогательную последовательность числовых функций.
\begin{dif}\label{dif:1.++}
Через \(\Lambda_m(\lambda)\), \(m=1,2,\ldots\), будут обозначаться функции
такие, что при любых фиксированных \(m\in\mathbb N\) и \(\lambda_0\in
(\sigma,\tau)\) число \(\Lambda_m(\lambda_0)\) представляет собой \(m\)-ое
снизу (с учетом кратности) собственное значение оператора \(F(\lambda_0)\).
\end{dif}
Определение~\ref{dif:1.++} корректно, поскольку операторы-значения
оператор-функции \(F(\lambda)\) ограничены снизу и имеют чисто дискретный
спектр.

Основное свойство функций \(\Lambda_m(\lambda)\) состоит в следующем.
\begin{lem}\label{lem:1.2.1}
При любом фиксированном \(m\in\mathbb N\) функция \(\Lambda_m(\lambda)\)
непрерывна на интервале \((\sigma,\tau)\).
\end{lem}
\begin{proof}
Очевидно, достаточно доказать, что при любом фиксированном
\(\lambda_0\in(\sigma,\tau)\) функция \(\Lambda_m(\lambda)\) непрерывна
на некоторой окрестности числа \(\lambda_0\).

Пусть вещественное число \(\mu\) и окрестность \(\mathcal U(\lambda_0)
\Subset(\sigma,\tau)\) числа \(\lambda_0\) таковы, что при любом фиксированном
\(\lambda_1\in\mathcal U(\lambda_0)\) оператор \(F(\lambda_1)-\mu\) является
положительным. Такие \(\mu\) и \(\mathcal U(\lambda_0)\) существуют по
определению оператор-функции \(F(\lambda)\). Обозначим через
\(\Theta_{m,\mu}(\lambda)\) числовую функцию окрестности \(\mathcal
U(\lambda_0)\), ставящую любому \(\lambda_1\in\mathcal U(\lambda_0)\) в
соответствие \(m\)-ое сверху (с учетом кратности) собственное значение
оператора \((F(\lambda_1)-\mu)^{-1}\).

Заметим теперь, что оператор-функция \((F(\lambda)-\mu)^{-1}\) непрерывна на
окрестности \(\mathcal U(\lambda_0)\) в смысле равномерной операторной
топологии "--- это с очевидностью следует из определения оператор-функции
\(F(\lambda)\). Поэтому числовая функция \(\Theta_{m,\mu}(\lambda)\)
непрерывна на окрестности \(\mathcal U(\lambda_0)\). Однако число \(\mu\)
выбрано таким образом, что справедливо тождество
\[
	\Lambda_m(\lambda)=\dfrac{1}{\Theta_{m,\mu}(\lambda)}+\mu.
\]
Поэтому функция \(\Lambda_m(\lambda)\) также непрерывна на окрестности
\(\mathcal U(\lambda_0)\). Тем самым лемма доказана.
\end{proof}
\begin{proof}[Доказательство теоремы~\ref{tm:1.44}.]
Непосредственно из определения функции \(\nu_F(\lambda)\) и функций
\(\Lambda_m(\lambda)\) следует, что не менее \(\nu_F(\xi_2)-\nu_F(\xi_1)\)
функций из последовательности \(\{\Lambda_m(\lambda)\}_{m=1}^{\infty}\)
принимают неотрицательное значение в точке \(\xi_1\) и отрицательное "--- в
точке \(\xi_2\). Из утверждаемой леммой~\ref{lem:1.2.1} непрерывности функций
\(\Lambda_m(\lambda)\) потому следует, что не менее \(\nu_F(\xi_2)-
\nu_F(\xi_1)\) функций из последовательности \(\{\Lambda_m(\lambda)
\}_{m=1}^{\infty}\) имеют нуль на полуинтервале \([\xi_1,\xi_2)\). Утверждение
теоремы теперь немедленно получается из определения функций
\(\Lambda_m(\lambda)\) и определения понятия собственного значения
оператор-функции \(F(\lambda)\).
\end{proof}

\subsection{Оценки числа собственных значений при выполнении условий
монотонности}
Если наложить на вид оператор-функции \(F(\lambda)\) дополнительные
ограничения, то появится возможность оценить величины \(\mathcal
N_F(\xi_1,\xi_2)\) не только снизу, но и сверху. При формулировке
соответствующих утверждений будет использоваться понятие \emph{замыкания
квадратичной формы}, определенное, например, 
в~\cite[глава~\(\mathrm{VI}\)]{Kato}.

\begin{tm}\label{tm:1.55}
Пусть определена форма-функция \(\mathfrak f(\lambda)[y]\) со следующими
свойствами.
\begin{itemize}
\item При любом фиксированном \(\lambda_0\in(\sigma,\tau)\) квадратичная форма
\(\mathfrak f(\lambda_0)[y]\) представляет собой замыкание определенной на
линейном подпространстве \(\mathfrak D(F(\lambda_0))\) квадратичной формы
\(\langle F(\lambda_0)y,y\rangle_{\mathfrak H}\).
\item Область определения квадратичной формы \(\mathfrak f(\lambda_0)[y]\) не
зависит от выбора числа \(\lambda_0\in(\sigma,\tau)\). Ниже эта не зависящая от
выбора \(\lambda_0\in(\sigma,\tau)\) область будет обозначаться просто через
\(\mathfrak D(\mathfrak f)\).
\item При любых фиксированных \(\lambda_1,\lambda_2\in(\sigma,\tau)\),
удовлетворяющих неравенству \(\lambda_1<\lambda_2\), определенная на линейном
подпространстве \(\mathfrak D(\mathfrak f)\) квадратичная форма
\[
	\mathfrak f(\lambda_1)[y]-\mathfrak f(\lambda_2)[y]
\]
является положительной.
\end{itemize}
Тогда для любого полуинтервала \([\xi_1,\xi_2)\Subset(\sigma,\tau)\)
имеет место равенство
\[
	\mathcal N_F(\xi_1,\xi_2)=\nu_F(\xi_2)-\nu_F(\xi_1).
\]
\end{tm}

Доказательство теоремы~\ref{tm:1.55} будет основано на следующем приводимом без
доказательства утверждении, представляющем собой простое следствие известного
вариационного принципа~\cite[глава~1, теорема~\(13^{bis}\)]{Gl}.
\begin{prop}\label{var.pr}
Пусть \(F\) есть ограниченный снизу самосопряженный оператор гильбертова
пространства \(\mathfrak H\), спектр которого чисто дискретен на луче
\((-\infty,0)\). Пусть \(\mathfrak f[y]\) есть замыкание определенной на
линейном подпространстве \(\mathfrak D(F)\) квадратичной формы
\(\langle Fy,y\rangle_{\mathfrak H}\). Тогда число отрицательных собственных
значений оператора \(F\) равно максимуму размерностей подпространств
\(\mathfrak M\subset\mathfrak D(\mathfrak f)\), на которых квадратичная форма
\(\mathfrak f[y]\) отрицательна.
\end{prop}

\begin{proof}[Доказательство теоремы~\ref{tm:1.55}.]
Зафиксируем произвольное число \(m\in\mathbb N\), а также произвольные
числа \(\lambda_1,\lambda_2\in(\sigma,\tau)\), удовлетворяющие неравенству
\(\lambda_1<\lambda_2\). Очевидно, что найдется \(m\)-мерное подпространство
\(\mathfrak M\subset\mathfrak D(F(\lambda_1))\), на котором квадратичная форма
оператора \(F(\lambda_1)-\Lambda_m(\lambda_1)\) будет неположительна. Однако
при выполнении условий теоремы указанное подпространство \(\mathfrak M\) будет
подпространством области определения квадратичной формы \(\mathfrak
f(\lambda_2)[y]\), причем квадратичная форма \(\mathfrak f(\lambda_2)[y]-
\Lambda_m(\lambda_1)\cdot\|y\|^2_{\mathfrak H}\) будет на подпространстве
\(\mathfrak M\) отрицательна. Из утверждения~\ref{var.pr} теперь немедленно
вытекает справедливость неравенства
\(\Lambda_m(\lambda_2)<\Lambda_m(\lambda_1)\).

Ввиду произвольности выбора чисел \(m\), \(\lambda_1\) и \(\lambda_2\),
сказанное означает, что все функции \(\Lambda_m(\lambda)\) строго убывают на
интервале \((\sigma,\tau)\).

Однако из непрерывности и строгого убывания функций \(\Lambda_m(\lambda)\)
следует, что величина \(\mathcal N_F(\xi_1,\xi_2)\) в точности равна числу
таких функций из последовательности \(\{\Lambda_m(\lambda)\}_{m=1}^{\infty}\),
которые принимают неотрицательное значение в точке \(\xi_1\) и отрицательное
"--- в точке \(\xi_2\). Поскольку число таких функций есть \(\nu_F(\xi_2)-
\nu_F(\xi_1)\), то утверждение теоремы справедливо.
\end{proof}

\subsection{Оценки числа собственных значений при выполнении условий
отрицательности типа спектра}
Условия теоремы~\ref{tm:1.55} можно в некотором смысле ослабить. Именно,
справедливо следующее утверждение.
\begin{tm}\label{tm:1.66}
Пусть определены формы-функции \(\mathfrak f(\lambda)[y]\) и \(\mathfrak
f'(\lambda)[y]\) со следующими свойствами.
\begin{itemize}
\item При любом фиксированном \(\lambda_0\in(\sigma,\tau)\) квадратичная форма
\(\mathfrak f(\lambda_0)[y]\) представляет собой замыкание определенной на
линейном подпространстве \(\mathfrak D(F(\lambda_0))\) квадратичной формы
\(\langle F(\lambda_0)y,y\rangle_{\mathfrak H}\).
\item Область определения квадратичной формы \(\mathfrak f(\lambda_0)[y]\) не
зависит от выбора числа \(\lambda_0\in(\sigma,\tau)\). Ниже эта не зависящая от
выбора \(\lambda_0\in(\sigma,\tau)\) область будет обозначаться просто через
\(\mathfrak D(\mathfrak f)\).
\item При любом фиксированном \(\lambda_0\in(\sigma,\tau)\) областью
определения квадратичной формы \(\mathfrak f'(\lambda_0)[y]\) является
линейное подпространство \(\mathfrak D(\mathfrak f)\).
\item Квадратичные формы \(\mathfrak f'(\lambda)[y]\) равномерно по
\(\lambda\in(\sigma,\tau)\) ограничены относительно произвольной метрики, по
которой \(\mathfrak D(\mathfrak f)\) полно.
\item При любом фиксированном \(y_0\in\mathfrak D(\mathfrak f)\) числовая
функция \(\mathfrak f(\lambda)[y_0]\) дифференцируема на интервале
\((\sigma,\tau)\), и производная этой функции есть числовая функция
\(\mathfrak f'(\lambda)[y_0]\).
\item Для любой собственной пары \(\{\lambda_0,y_0\}\) оператор-функции
\(F(\lambda)\) справедливо неравенство \(\mathfrak f'(\lambda_0)[y_0]<0\).
\end{itemize}
Тогда для любого полуинтервала \([\xi_1,\xi_2)\Subset(\sigma,\tau)\)
справедливо равенство
\[
	\mathcal N_F(\xi_1,\xi_2)=\nu_F(\xi_2)-\nu_F(\xi_1).
\]
\end{tm}

Отметим, что содержание теоремы~\ref{tm:1.66} тесно связано с известным в
теории операторных пучков понятием \emph{типа собственного значения}.

\begin{proof}[Доказательство теоремы~\ref{tm:1.66}.]
Пусть \(m\in\mathbb N\) и \(\lambda_0\in(\sigma,\tau)\) "--- произвольные
числа, для которых справедливо равенство \(\Lambda_m(\lambda_0)=0\).

Введем в линейном пространстве \(\mathfrak D(\mathfrak f)\) норму
\[
	\|y\|^2_{\mathfrak F}:=\mathfrak f(\lambda_0)[y]-\mu
	\cdot\|y\|^2_{\mathfrak H},
\]
где \(\mu\) есть произвольно фиксированное вещественное число, строго меньшее
вершины формы \(\mathfrak f(\lambda_0)[y]\). Относительно нормы
\(\|y\|_{\mathfrak F}\) линейное подпространство \(\mathfrak D(\mathfrak f)\)
является гильбертовым пространством, которое ниже будет обозначаться через
\(\mathfrak F\).

Рассмотрим оператор-функцию \(K(\lambda)\), для которой при любом фиксированном
\(\lambda_1\in(\sigma,\tau)\) оператор \(K(\lambda_1)\) действует в
гильбертовом пространстве \(\mathfrak F\) и удовлетворяет тождеству
\[
	\langle (1+K(\lambda_1))y,y\rangle_{\mathfrak F}=
	\mathfrak f(\lambda_1)[y].
\]
Из теоремы о замкнутом графике (см., например,~\cite[глава~\(\mathrm{III}\),
теорема~\mbox{5.20}]{Kato}) нетрудно вывести тот факт, что значениями
оператор-функции \(K(\lambda)\) являются ограниченные самосопряженные
операторы. Кроме того, из определения нормы \(\|y\|_{\mathfrak F}\) и
компактности резольвенты оператора \(F(\lambda_0)\) вытекает факт компактности
оператора \(K(\lambda_0)\).

Заметим, что оператор-функция \(K(\lambda)\) является дифференцируемой, причем
ее производная \(K'(\lambda)\) удовлетворяет тождеству
\[
	\langle K'(\lambda)y,y\rangle_{\mathfrak F}=\mathfrak f'(\lambda)[y].
\]
Заметим также, что операторы-значения оператор-функции \(K'(\lambda)\)
равномерно по \(\lambda\in(\sigma,\tau)\) ограничены в силу условий теоремы.

Разложим теперь пространство \(\mathfrak F\) в прямую сумму \(\mathfrak
F_-\oplus\mathfrak F_0\oplus\mathfrak F_+\), где
\begin{itemize}
\item \(\mathfrak F_-\) есть инвариантное подпространство оператора
\(K(\lambda_0)\), отвечающее той части его спектра, которая лежит на луче
\((-\infty,-1)\);
\item \(\mathfrak F_0\) есть инвариантное подпространство оператора
\(K(\lambda_0)\), отвечающее собственному значению \(-1\);
\item \(\mathfrak F_+\) есть инвариантное подпространство оператора
\(K(\lambda_0)\), отвечающее той части его спектра, которая лежит на луче
\((-1,+\infty)\).
\end{itemize}
Из компактности оператора \(K(\lambda_0)\) следует существование такого
\(\varepsilon>0\), для которого на подпространстве \(\mathfrak F_-\) будет
справедлива оценка
\begin{equation}\label{eq:1.616}
	\langle K(\lambda_0)y,y\rangle_{\mathfrak F}\leqslant
	(-1-\varepsilon)\|y\|^2_{\mathfrak F},
\end{equation}
а на подпространстве \(\mathfrak F_+\) будет справедлива оценка
\begin{equation}\label{eq:1.666}
	\langle K(\lambda_0)y,y\rangle_{\mathfrak F}\geqslant
	(-1+\varepsilon)\|y\|^2_{\mathfrak F}.
\end{equation}

Однако из условий теоремы и из первой теоремы о представлении
(см.~\cite[глава~\(\mathrm{VI}\), теорема~\mbox{2.1}]{Kato}) следует, что для
любого ненулевого вектора \(y\in\mathfrak F_0\) выполняется неравенство 
\[
	\langle K'(\lambda_0)y,y\rangle_{\mathfrak F}<0.
\]
Поэтому оценка~\eqref{eq:1.616} означает, что найдется правая проколотая
полуокрестность \(\mathcal U_{\mbox{\footnotesize пр}}(\lambda_0)
\Subset(\sigma,\tau)\) точки \(\lambda_0\) такая, что для любого \(\lambda_1
\in\mathcal U_{\mbox{\footnotesize пр}}(\lambda_0)\) квадратичная форма
оператора \(1+K(\lambda_1)\) будет отрицательна на подпространстве \(\mathfrak
F_-\oplus\mathfrak F_0\). Аналогично, из оценки~\eqref{eq:1.666} и равномерной
ограниченности оператор-функции \(K'(\lambda)\) следует, что найдется левая
проколотая полуокрестность \(\mathcal U_{\mbox{\footnotesize лев}}(\lambda_0)
\Subset(\sigma,\tau)\) точки \(\lambda_0\) такая, что для любого \(\lambda_1\in
\mathcal U_{\mbox{\footnotesize лев}}(\lambda_0)\) квадратичная форма оператора
\(1+K(\lambda_1)\) будет положительна на подпространстве \(\mathfrak F_0\oplus
\mathfrak F_+\).

Заметим теперь, что подпространство \(\mathfrak F_-\oplus\mathfrak F_0\) имеет
размерность \(\geqslant m\). Поэтому из определения оператор-функции
\(K(\lambda)\) и утверждения~\ref{var.pr} следует, что при любом \(\lambda_1
\in\mathcal U_{\mbox{\footnotesize пр}}(\lambda_0)\) справедливо неравенство
\(\Lambda_m(\lambda_1)<0\). Аналогично, подпространство \(\mathfrak
F_0\oplus\mathfrak F_+\) имеет коразмерность \(<m\). Поэтому при любом
\(\lambda_1\in\mathcal U_{\mbox{\footnotesize лев}}(\lambda_0)\) справедливо
неравенство \(\Lambda_m(\lambda_1)\geqslant 0\).

Итак, получен следующий результат: \emph{каковы бы ни были \(m\in\mathbb N\) и
\(\lambda_0\in(\sigma,\tau)\), удовлетворяющие равенству
\(\Lambda_m(\lambda_0)=0\), функция \(\Lambda_m(\lambda)\) неотрицательна в
некоторой левой полуокрестности точки \(\lambda_0\) и отрицательна в некоторой
проколотой правой полуокрестности точки \(\lambda_0\).} Поскольку любая из
функций \(\Lambda_m(\lambda)\) непрерывна в силу леммы~\ref{lem:1.2.1}, то
этот результат означает, что любая из функций \(\Lambda_m(\lambda)\) строго
убывает в каждом из своих нулей. Отсюда, в свою очередь, следует, что каждая
из функций \(\Lambda_m(\lambda)\) имеет не более одного нуля, слева от
которого она принимает только положительные, а справа "--- только
отрицательные значения.

Повторяя теперь рассуждения, завершающие доказательство теоремы~\ref{tm:1.55},
убеждаемся в справедливости утверждения теоремы.
\end{proof}

\section{Применение к дифференциальным оператор-функциям}
Рассмотрим теперь оператор-функции, описываемые следующим определением.
\begin{dif}\label{dif:1.00}
Пусть \(p_k(x,\lambda)\), \(k=0,\ldots,n\), есть непрерывные вещественные
функции прямоугольника \([a,b]\times(\sigma,\tau)\). Пусть при этом для любого
фиксированного \(\lambda_0\in(\sigma,\tau)\) функция \(p_0(x,\lambda_0)\)
положительна и принадлежит пространству \(C^n[a,b]\), а любая из функций
\(p_k(x,\lambda_0)\), \(k=1,\ldots,n-1\), принадлежит пространству
\(C^{n-k}[a,b]\). Пусть, далее, \(U(\lambda)\) есть непрерывная
матрица-функция интервала \((\sigma,\tau)\), значениями которой являются
унитарные комплексные матрицы размера \(2n\times 2n\). Тогда
оператор-функцией, определенной дифференциальным выражением
\begin{equation}\label{eq:1.1.1}
	\sum\limits_{k=0}^n(-1)^{n-k}\left(p_k(x,\lambda)
	y^{(n-k)}(x)\right)^{(n-k)}
\end{equation}
и краевыми условиями
\begin{equation}\label{eq:1.1.2}
	(U(\lambda)-1)y^{\vee}+i(U(\lambda)+1)y^{\wedge}=0,
\end{equation}
называется оператор-функция \(S(\lambda)\) со следующими свойствами.
\begin{itemize}
\item При любом фиксированном \(\lambda_0\in(\sigma,\tau)\) оператор
\(S(\lambda_0)\) действует в гильбертовом пространстве \(L_2[a,b]\).
\item При любом фиксированном \(\lambda_0\in(\sigma,\tau)\) областью
определения оператора \(S(\lambda_0)\) является плотное в гильбертовом
пространстве \(L_2[a,b]\) линейное подпространство
\[
	\left\{y\in L_2[a,b]\;\vline\;y\in W_2^{2n}[a,b],\,
	(U(\lambda_0)-1)y^{\vee}+i(U(\lambda_0)+1)y^{\wedge}=0\right\},
\]
где \(y^{\wedge}\) и \(y^{\vee}\) есть связанные с функцией \(y(x)\) векторы
пространства \(\mathbb C^{2n}\) вида
\begin{align*}
	y^{\wedge}&:=\begin{pmatrix} y(a)\\ y'(a)\\ \hdotsfor{1}
	\\ y^{(n-1)}(a)\\ y(b)\\ y'(b)\\ \hdotsfor{1}\\
	y^{(n-1)}(b)\end{pmatrix},&
	y^{\vee}&:=\begin{pmatrix} y^{[2n-1]}(a)\\
	y^{[2n-2]}(a)\\ \hdotsfor{1}\\ y^{[n]}(a)\\
	-y^{[2n-1]}(b)\\ -y^{[2n-2]}(b)\\ \hdotsfor{1}\\
	-y^{[n]}(b)\end{pmatrix}
\end{align*}
(ср.~с~\cite[(7.50)]{Rofe}). Здесь через \(y^{[n+m]}(x)\) обозначены хорошо
известные в теории дифференциальных операторов квазипроизводные вида
\begin{align*}
	y^{[n+m]}(x)&:=\sum\limits_{k=0}^m(-1)^{m-k}
	\left(p_k(x,\lambda_0)y^{(n-k)}(x)\right)^{(m-k)},&
	m&=0,\ldots,n
\end{align*}
(ср.~с~\cite[\S\,15, (3)]{Najm} и~\cite[(7.46)]{Rofe}).
\item При любом фиксированном \(\lambda_0\in(\sigma,\tau)\) действие оператора
\(S(\lambda_0)\) на произвольную функцию \(y\in\mathfrak D(S(\lambda_0))\)
определяется равенством
\[
	S(\lambda_0)y:=\sum\limits_{k=0}^n(-1)^{n-k}\left(
	p_k(x,\lambda_0)y^{(n-k)}(x)\right)^{(n-k)}.
\]
\end{itemize}
\end{dif}
Настоящий пункт посвящается применению теорем~\ref{tm:1.44}\,--\,\ref{tm:1.66}
к оператор-функциям \(S(\lambda)\), описываемым определением~\ref{dif:1.00}.

\subsection{Свойства дифференциальных оператор-функций}
Оператор-функции, описываемые определением~\ref{dif:1.00}, обладают следующими
свойствами.
\begin{lem}\label{lem:1.1}
Пусть \(A(\lambda)\) есть матрица-функция, значение которой при любом
\(\lambda_0\in(\sigma,\tau)\) удовлетворяет равенству
\begin{equation}\label{eq:1.33.3}
	A(\lambda_0)=\dfrac{1}{2\pi}\int\limits_{\Gamma}
	\dfrac{z+1}{z-1}\cdot(U(\lambda_0)-z)^{-1}\,dz.
\end{equation}
Здесь \(U(\lambda_0)\) есть унитарная матрица из определяющих оператор
\(S(\lambda_0)\) краевых условий, а \(\Gamma\) есть положительно
ориентированный контур, внутренность которого не содержит точку \(1\), но
содержит все отличные от \(1\) собственные значения матрицы \(U(\lambda_0)\).
Тогда при любом фиксированном \(\lambda_0\in(\sigma,\tau)\) для любой функции
\(y\in\mathfrak D(S(\lambda_0))\) справедливо равенство
\begin{equation}\label{eq:1.1.4}
	\langle S(\lambda_0)y,y\rangle_{L_2[a,b]}=\sum\limits_{k=0}^n
	\int\limits_a^b p_k(x,\lambda_0)\left|y^{(n-k)}(x)\right|^2\,
	dx+\langle A(\lambda_0)y^{\wedge},y^{\wedge}\rangle_{\mathbb C^{2n}}.
\end{equation}
\end{lem}

\begin{lem}\label{tm:1.00}
Пусть матрица-функция \(U(\lambda)\) такова, что числовая функция
\(\operatorname{rank}(U(\lambda)-1)\) постоянна на интервале
\((\sigma,\tau)\). Тогда для любого фиксированного
\(\lambda_0\in(\sigma,\tau)\) найдутся такое вещественное число \(\mu\) и
такая окрестность \(\mathcal U(\lambda_0)\Subset(\sigma,\tau)\) числа
\(\lambda_0\), что при любом фиксированном \(\lambda_1\in\mathcal
U(\lambda_0)\) оператор \(S(\lambda_1)-\mu\) будет неотрицателен.
\end{lem}

\begin{lem}\label{lem:1.21}
Оператор-функция \((S(\lambda)-i)^{-1}\) непрерывна на интервале
\((\sigma,\tau)\) в смысле равномерной операторной топологии.
\end{lem}

\begin{lem}\label{lem:1.33}
При любом фиксированном \(\lambda_0\in(\sigma,\tau)\) замыканием определенной
на линейном подпространстве \(\mathfrak D(S(\lambda_0))\) квадратичной формы
\(\langle S(\lambda_0)y,y\rangle_{L_2[a,b]}\) является квадратичная форма
\(\mathfrak s(\lambda_0)[y]\), областью определения которой является линейное
подпространство
\[
	\left\{y\in L_2[a,b]\;\vline\;y\in W_2^n[a,b],\,y^{\wedge}\perp
	\ker(U(\lambda_0)-1)\,\right\},
\]
и значение которой на произвольной функции \(y\in\mathfrak D(\mathfrak s)\)
определяется равенством
\[
	\mathfrak s(\lambda_0)[y]:=\sum\limits_{k=0}^n
	\int\limits_a^b p_k(x,\lambda_0)\left|y^{(n-k)}(x)\right|^2\,
	dx+\langle A(\lambda_0)y^{\wedge},y^{\wedge}\rangle_{\mathbb C^{2n}}.
\]
\end{lem}

\begin{proof}[Доказательство леммы~\ref{lem:1.1}.]
При помощи интегрирования по частям нетрудно установить справедливость
тождества
\[
	\langle S(\lambda_0)y,y\rangle_{L_2[a,b]}=\sum\limits_{k=0}^n
	\int\limits_a^b p_k(x,\lambda_0)\left|y^{(n-k)}(x)\right|^2\,
	dx+\langle y^{\vee},y^{\wedge}\rangle_{\mathbb C^{2n}}.
\]
Поэтому лемма будет доказана, если будет установлено, что для любых векторов
\(Y,Z\in\mathbb C^{2n}\), удовлетворяющих условию
\begin{equation}\label{eq:1.33.5}
	(U(\lambda_0)-1)Z+i(U(\lambda_0)+1)Y=0,
\end{equation}
выполняется равенство
\begin{equation}\label{eq:1.33.6}
	\langle Z,Y\rangle_{\mathbb C^{2n}}=\langle A(\lambda_0)Y,
	Y\rangle_{\mathbb C^{2n}}.
\end{equation}

Заметим, что условие~\eqref{eq:1.33.5} можно переписать в параметрической форме
\begin{equation}\label{eq:1.33.55}
\begin{aligned}
	Y&=(U(\lambda_0)-1)X,\\
	Z&=-i(U(\lambda_0)+1)X,\\
	X&\in\mathbb C^{2n}.
\end{aligned}
\end{equation}
При этом левая часть предполагаемого равенства~\eqref{eq:1.33.6} оказывается
равна скалярному произведению
\[
	\langle B(\lambda_0)X,X\rangle_{\mathbb C^{2n}},
\]
где через \(B(\lambda_0)\) обозначена матрица
\[
	B(\lambda_0):=-i\left(U^{-1}(\lambda_0)-1\right)\cdot
	(U(\lambda_0)+1)=\\=\dfrac{1}{2\pi}\int\limits_{\Gamma_1}
	(z^{-1}-1)\cdot(z+1)\cdot(U(\lambda_0)-z)^{-1}\,dz.
\]
Здесь \(\Gamma_1\) есть положительно ориентированная граница кольца
\(\{z\in\mathbb C\;\vline\;|z|\in(1-\varepsilon,1+\varepsilon)\}\).
С другой стороны, правая часть предполагаемого равенства~\eqref{eq:1.33.6}
оказывается равна скалярному произведению
\[
	\langle C(\lambda_0)X,X\rangle_{\mathbb C^{2n}},
\]
где через \(C(\lambda_0)\) обозначена матрица
\[
	C(\lambda_0):=\dfrac{1}{2\pi}\int\limits_{\Gamma}(z^{-1}-1)
	\cdot(z+1)\cdot(U(\lambda_0)-z)^{-1}\,dz.
\]
Здесь \(\Gamma\) есть контур из определения матрицы \(A(\lambda_0)\).

Между тем, функция \(-i(z^{-1}-1)(z+1)\) обращается в нуль в точке \(z=1\), а
матрица \(U(\lambda_0)\) заведомо не имеет присоединенных векторов, отвечающих
собственному значению \(1\). Поэтому матрицы \(B(\lambda_0)\) и
\(C(\lambda_0)\) равны. Следовательно, при выполнении
условия~\eqref{eq:1.33.5} равенство~\eqref{eq:1.33.6} действительно
справедливо. Тем самым лемма доказана.
\end{proof}

\begin{proof}[Доказательство леммы~\ref{tm:1.00}.]
Зафиксируем вещественное число \(\mu\) таким образом, что определенная в
пространстве \(W_2^n[a,b]\) ограниченная квадратичная форма
\[
	\mathfrak s_{\lambda_0,\mu}[y]:=\sum\limits_{k=0}^n
	\int\limits_a^b p_k(x,\lambda_0)\left|y^{(n-k)}(x)\right|^2\,
	dx-\int\limits_a^b \mu\cdot|y(x)|^2\,dx+\\+
	\langle A(\lambda_0)y^{\wedge},y^{\wedge}\rangle_{\mathbb C^{2n}}
\]
является равномерно положительной.

Заметим, что при выполнении условия \(\operatorname{rank}(U(\lambda)-1)=
\mathrm{const}\) матрица-функция \(A(\lambda)\) является непрерывной на
интервале \((\sigma,\tau)\). Поэтому для чисел \(\lambda_1\), достаточно
близких к \(\lambda_0\), разность между
квадратичной формой \(\mathfrak s_{\lambda_0,\mu}[y]\) и квадратичной формой
\[
	\mathfrak s_{\lambda_1,\mu}[y]:=\sum\limits_{k=0}^n
	\int\limits_a^b p_k(x,\lambda_1)\left|y^{(n-k)}(x)\right|^2\,
	dx-\int\limits_a^b \mu\cdot|y(x)|^2\,dx+\\+
	\langle A(\lambda_1)y^{\wedge},y^{\wedge}\rangle_{\mathbb C^{2n}}
\]
будет равномерно мала (по метрике пространства \(W_2^n[a,b]\)). Следовательно,
найдется такая окрестность \(\mathcal U(\lambda_0)\Subset(\sigma,\tau)\) числа
\(\lambda_0\), что для любого фиксированного \(\lambda_1\in\mathcal
U(\lambda_0)\) квадратичная форма \(\mathfrak s_{\lambda_1,\mu}[y]\) будет
неотрицательна.

Однако из леммы~\ref{lem:1.1} с очевидностью следует, что неотрицательность
квадратичной формы \(\mathfrak s_{\lambda_1,\mu}[y]\) влечет за собой
неотрицательность оператора \(S(\lambda_1)-\mu\) (обратное, вообще говоря,
неверно). Значит, ранее определенные число \(\mu\) и окрестность \(\mathcal
U(\lambda_0)\) обладают тем свойством, что для любого фиксированного
\(\lambda_1\in\mathcal U(\lambda_0)\) оператор \(S(\lambda_1)-\mu\) является
неотрицательным. Тем самым утверждение леммы справедливо.
\end{proof}

\begin{proof}[Доказательство леммы~\ref{lem:1.21}.]
Рассмотрим определенную на параллелепипеде \([a,b]\times[a,b]\times
(\sigma,\tau)\) функцию \(G(x,t,\lambda)\), для которой при любом
фиксированном \(\lambda_0\in(\sigma,\tau)\) функция \(G(x,t,\lambda_0)\)
представляет собой функцию Грина оператора \(S(\lambda_0)-i\). Нетрудно
проверить, что непрерывность функций \(p_k(x,\lambda)\) и \(U(\lambda)\)
приводит к тому, что функция \(G(x,t,\lambda)\) непрерывна на параллелепипеде
\([a,b]\times[a,b]\times(\sigma,\tau)\). Отсюда с очевидностью вытекает
справедливость утверждения леммы.
\end{proof}

\begin{proof}[Доказательство леммы~\ref{lem:1.33}.]
Справедливость утверждения доказываемой леммы нетрудно вывести из
леммы~\ref{lem:1.1} и возможности записи условия~\eqref{eq:1.33.5} в
форме~\eqref{eq:1.33.55}.
\end{proof}

\subsection{Оценки числа собственных значений}
Из лемм~\ref{lem:1.1}\,--\,\ref{lem:1.33} и
теорем~\ref{tm:1.44}\,--\,\ref{tm:1.66} вытекает справедливость следующих
утверждений.

\begin{tm}\label{tm:1.4}
Пусть матрица-функция \(U(\lambda)\) такова, что числовая функция
\(\operatorname{rank}(U(\lambda)-1)\) постоянна на интервале
\((\sigma,\tau)\). Тогда для любого полуинтервала \([\xi_1,\xi_2)\Subset
(\sigma,\tau)\) имеет место оценка
\begin{equation}\label{eq:1.12+}
	\mathcal N_S(\xi_1,\xi_2)\geqslant\nu_S(\xi_2)-\nu_S(\xi_1).
\end{equation}
\end{tm}
\begin{tm}\label{tm:1.5}
Пусть выполнены следующие условия.
\begin{itemize}
\item Подпространство \(\ker\left(U(\lambda_0)-1\right)\) не зависит от выбора
числа \(\lambda_0\in(\sigma,\tau)\).
\item При любом фиксированном \(x\in[a,b]\) функции \(p_k(x,\lambda)\),
\(k=0,\ldots,n-1\), невозрастают по \(\lambda\). Кроме того, при любом
фиксированном \(x\in[a,b]\) функция \(p_n(x,\lambda)\) строго убывает по
\(\lambda\).
\item Отвечающая матрице-функции \(U(\lambda)\) матрица-функция \(A(\lambda)\)
(см. лемму~\ref{lem:1.1}) невозрастает по \(\lambda\) в том смысле, что при
любых фиксированных \(\lambda_1,\lambda_2\in(\sigma,\tau)\), удовлетворяющих
неравенству \(\lambda_1<\lambda_2\), самосопряженная матрица \(A(\lambda_1)-
A(\lambda_2)\) является неотрицательной.
\end{itemize}
Тогда для любого полуинтервала \([\xi_1,\xi_2)\Subset(\sigma,\tau)\)
имеет место равенство
\[
	\mathcal N_S(\xi_1,\xi_2)=\nu_S(\xi_2)-\nu_S(\xi_1).
\]
\end{tm}
\begin{tm}\label{tm:1.6}
Пусть выполнены следующие условия.
\begin{itemize}
\item Подпространство \(\ker(U(\lambda_0)-1)\) не зависит от выбора числа
\(\lambda_0\in(\sigma,\tau)\).
\item Функции \(p_k(x,\lambda)\), \(k=0,\ldots,n\), являются дифференцируемыми
по параметру \(\lambda\), причем все производные
\(\{p_k(x,\lambda)\}'_{\lambda}\) непрерывны на прямоугольнике
\([a,b]\times(\sigma,\tau)\).
\item Матрица-функция \(U(\lambda)\) является непрерывно дифференцируемой.
\item Для любой собственной пары \(\{\lambda_0,y_0\}\) оператор-функции
\(S(\lambda)\) имеет место неравенство
\[
	\sum\limits_{k=0}^n\int\limits_a^b\left.
	\left\{p_k(x,\lambda)\right\}'_{\lambda}\right|_{\lambda=\lambda_0}
	\cdot\left|y^{(n-k)}(x)\right|\,dx+\langle A'(\lambda_0)y^{\wedge},
	y^{\wedge}\rangle_{\mathbb C^{2n}}<0.
\]
\end{itemize}
Тогда для любого полуинтервала \([\xi_1,\xi_2)\Subset(\sigma,\tau)\)
имеет место равенство
\[
	\mathcal N_S(\xi_1,\xi_2)=\nu_S(\xi_2)-\nu_S(\xi_1).
\]
\end{tm}

В качестве следствий теорем~\ref{tm:1.4}\,--\,\ref{tm:1.6} могут быть получены
некоторые известные оценки числа собственных значений дифференциальных
оператор-функций "--- в
частности, оценки, связанные с результатами работ~\cite{Kamke} и~\cite{MSS}.

Автор выражает признательность проф.~А.~А.~Шкаликову за постановку задач и
ценные замечания.

\end{document}